\documentclass[12pt]{article}
\usepackage{geometry}
\usepackage{chngpage}
\usepackage{graphicx}
\usepackage{amsmath}
\usepackage{amsfonts}
\usepackage{amssymb}
\usepackage{latexsym}
\usepackage[brazil]{varioref}
\usepackage[english,brazil]{babel}
\usepackage[latin1]{inputenc}
\usepackage[dvips]{color}

\newtheorem{defn}{Defini\c{c}{\~a}o}[subsection]
\newtheorem{teo}{Teorema}[subsection]
\newtheorem{coro}[teo]{Corol{\'a}rio}
\newtheorem{lema}[teo]{Lema}
\newtheorem{prop}[teo]{Proposição}
\newenvironment{dem}[1][Demonstra\c{c}{\~a}o]{\textbf{#1:}\ }
{\hfill\rule{1ex}{1ex}}

\begin{document}

\title{Um modelo algébrico do quantificador da ubiquidade}
\author{Tiago Augusto dos Santos Boza  \thanks{ Email: boza.tiago@gmail.com. Pós-Graduação em Filosofia, UNESP, FFC - Marília} \\
Hércules de Araujo Feitosa \thanks{Email: haf@fc.unesp.br. Departamento de Matemática, UNESP, FC - Bauru}}

\maketitle \abstract{O objetivo desse artigo é um aprofundamento das noções de quantificação
dentro do aspecto das lógicas moduladas. Para tanto, aborda-se a lógica
modulada do plausível, que procura formalizar o quantificador da ubiquidade. O
texto apresenta uma proposta, introduzida por Paul Halmos, de interpretação
da lógica quantificacional clássica em modelos algébricos e, como contribuição
original, estende este modelo para um modelo algébrico para a lógica do
plausível.

\medskip

\noindent{\bf Palavras Chave:} quantificadores; lógicas moduladas; lógica algébrica;
ubiquidade.}


\section*{Introdução}
${}$ \hspace{0,5 cm} Este artigo pretende apresentar algumas reflexões sobre
quantificadores, clássicos e não clássicos, suas formalizações e interpretações semânticas.

Para tanto, em particular, trataremos de um quantificador modulado, não clássico
que, portanto, não pode ser definido a partir dos quantificadores clássicos, no nosso caso, o quantificador da ubiquidade, que será apresentado na Seção 1, com sua respectiva interpretação semântica.

Na Seção seguinte, apresentamos uma forma distinta de interpretação para os quantificadores, como um operador numa álgebra de Boole. Este conceito será apresentado num ambiente algébrico, a partir das funções proposicionais e, por fim, das lógicas e álgebras monádicas. Nos últimos resultados desta Seção, apresentamos resultados que garantem que tal versão é adequada à lógica de primeira ordem e, desse modo,
equivalente aos modelos usuais para a lógica clássica de primeira ordem.

Como elemento original deste artigo, apresentamos, na última Seção, uma interpretação algébrica do quantificador da ubiquidade que estende as álgebras monádicas da Seção anterior.


\section{Lógica do Plausível}
${}$ \hspace {0,5 cm} Nesse capítulo, apresentaremos uma lógica modulada, a saber, a lógica do plausível.


\subsection{Versão Axiomática}

${}$ \hspace {0,5 cm} Uma lógica modulada é uma extensão conservativa da lógica clássica de primeira ordem, dada pela inclusão de um novo quantificador, não definível a partir dos quantificadores clássicos, universal $\forall$ e existencial $\exists$.

Para a lógica do plausível, inclui-se um novo quantificador na linguagem da lógica clássica de primeira ordem, denotado por $U$, denominado de quantificador do plausível ou quantificador da ubiquidade, que assim o chamaremos neste artigo.

Este quantificador tem a motivação intuitiva de capturar proposições da linguagem natural da forma "uma boa parte". Para isso, $\cite{gec08}$ introduziu como uma interpretação deste quantificador uma estrutura matemática denominada de "espaço pseudotopológico, a qual será detalhada numa próxima seção deste capítulo.

Desse modo, uma sentença quantificada como $UxAx$ significa que "uma boa parte de $x$, satisfaz a propriedade $A$" ou também que "há suficientes $x$ tais que $Ax$".

Agora formalmente, consideremos então $\mathcal{L}$ a lógica clássica de primeira ordem com identidade,
como em $\cite{fep05}$. A lógica do plausível $\mathcal{L}(U)$ é determinada a partir de $\mathcal{L}$ do seguinte modo, como em $\cite{gec08}$:

\begin{defn}
$\mathcal{L}(U)$ é determinada por todos os axiomas de $\mathcal{L}$ acrescidos dos seguintes axiomas específicos para o novo quantificador U: \\

$(Ax_1) (UxAx \land UxBx) \rightarrow Ux(Ax \land Bx)$ \\

$(Ax_2)$ $(UxAx \land UxBx) \rightarrow Ux(Ax \lor Bx)$ \\

$(Ax_3)$ $\forall xAx \rightarrow UxAx$ \\

$(Ax_4)$ $UxAx \rightarrow \exists xAx$ \\

$(Ax_5)$ $(\forall x(Ax \leftrightarrow Bx)) \rightarrow (UxAx \leftrightarrow UxBx)$ \\

$(Ax_6)$ $UxAx \rightarrow UyAy$, quando y ocorre livre para x em A. \\

- As regras de dedução são as seguintes: \\

Modus Ponens (MP): A, $A \rightarrow B$ $\vdash$ B \\

Generalização (Gen): $A \vdash \forall x A$
\end{defn}

As definições sintáticas usuais para $\mathcal{L}(U)$ como sentença, demonstração, teorema, consistência, entre outras, são definidas da mesma forma que na lógica clássica de primeira ordem.

Abaixo apenas apresentaremos alguns resultados desta lógica, cujas demonstrações podem ser encontradas em $\cite{gec08}$.

\begin{teo}
As fórmulas abaixo são teoremas de $\mathcal{L}(U)$: \\

(i) $Ux (Ax \lor \neg Ax)$; \\

(ii) $UxAx \land UxBx \rightarrow \exists x(Ax \land Bx)$; \\

(iii) $UxAx \rightarrow \neg Ux \neg Ax$.
\end{teo}

\begin{teo}
O cálculo de predicados $\mathcal{L}(U)$ é consistente. \\
\begin{dem}
Análoga à demonstração para a lógica clássica $\mathcal{L}$, mudando apenas a definição da função esquecimento, h, por meio da inclusão da condição h(UxAx) = h(Ax).
\end{dem}
\end{teo}

\begin{teo} (Teorema da Dedução) Seja $\Gamma \cup \lbrace A, B \rbrace$ um conjunto de
fórmulas de $\mathcal{L}(U)$. Suponhamos que $\Gamma \cup \lbrace A \rbrace \vdash B$, e que $x_i$ é uma variável livre de A e que, na demonstração de B a partir de $\Gamma \cup \lbrace A \rbrace$, a regra (Gen) não é aplicada em nenhuma fórmula $A_i$, que dependa de A. Neste caso, $\Gamma \vdash A \rightarrow B$.
\end{teo}

\begin{teo}
Seja $\Gamma$ um conjunto de fórmulas de $\mathcal{L}(U)$. Então $\Gamma$ é consistente see todo subconjunto finito $\Gamma_0$ de $\Gamma$ é consistente.
\end{teo}

\begin{teo}
Sejam $\Gamma$ um conjunto de fórmulas e A uma sentença de $\mathcal{L}(U)$. Então $\Gamma \cup \lbrace A \rbrace$ é inconsistente see $\Gamma \vdash \neg A$.
\end{teo}

\begin{teo}
Se $\Gamma$ é um conjunto de fórmulas de $\mathcal{L}(U)$ consistente maximal e A, B são sentenças de $\mathcal{L}(U)$, então: \\

(i) $\Gamma \vdash A$ see $A \in \Gamma$; \\

(ii) $A \notin \Gamma$ see $\neg A \in \Gamma$; \\

(iii) $A \land B \in \Gamma$ see $A, B \in \Gamma$.
\end{teo}

\begin{teo}
Todo conjunto consistente de sentenças de $\mathcal{L}(U)$ está contido em um conjunto consistente maximal.
\end{teo}

O resultado acima é similar ao Teorema de Lindenbaum.

Apresentaremos, na seção seguinte, a estrutura utilizada para interpretar tais elementos.


\subsection{Espaços Pseudo-Topológicos}

${}$ \hspace{0,5 cm} Nesta seção apresentamos os espaços pseudo-topológicos, que são ambientes matemáticos para a interpretação do quantificador da ubiquidade. $\cite{gec08}$, na tentativa de formalizar o novo quantificador $U$, procuraram uma estrutura matemática que pudesse modelar o seu novo quantificador. 

A ideia para esta nova estrutura matemática, obviamente, é advinda dos espaços topológicos. Há de se tomar cuidado com este fato, pois, apesar da semelhança conceitual, são estruturas distintas.

\begin{defn}
Um espaço pseudo-topológico é um par $(E, \Omega)$, em que E é um conjunto não-vazio e $\Omega$ um subconjunto de $\mathfrak{P}(E)$. Cada membro de $\Omega$ é denominado um aberto de $(E, \Omega)$, de maneira que: \\

(i) se $A, B \in \Omega$, então $A \cap B \in \Omega$; \\

(ii) se $A, B \in \Omega$, então $A \cup B \in \Omega$; \\

(iii) $E \in \Omega$; \\

(iv) $\emptyset \notin \Omega$. \\

Um subconjunto F de E é fechado em $(E, \Omega)$ quando seu complementar é um aberto em $(E, \Omega)$, isto é, o complementar de F, denotado por $F^C$, pertence à $\Omega$, $F^C \in \Omega$.
\end{defn}

Usualmente, uma topologia caracteriza uma operação de interior. Não tratamos com qualquer topologia, mas com pseudo-topologias. Assim, temos uma particular definição de interior.

A definição abaixo é apresentada originalmente neste trabalho, a fim de que as relações algébricas apresentadas posteriormente possuam uma interpretação nos espaços pseudo-topológicos.

\begin{defn}
Seja $(E, \Omega)$ um espaço pseudo-topológico. O interior de $A \subseteq \Omega$ é denotado por $\Upsilon(A)$ e definido por: \\

- o maior $D \in \Omega$ : $D \subseteq A$, se existe tal D; \\

- $\emptyset$, se não existe D.
\end{defn}

Segue desta definição que $\Upsilon(A) \subseteq A$, para todo $A \in \mathfrak{P}(E)$.

\begin{defn}
Seja $(E, \Omega)$ um espaço pseudo-topológico. Um conjunto $A \subseteq E$ é aberto em $(E, \Omega)$ se $A \subseteq \Upsilon(A)$.
\end{defn}

Se A é aberto então $\Upsilon(A) = A$. Além disso, o interior de um conjunto é um conjunto aberto ou é o conjunto vazio.

Abaixo, alguns resultados específicos para o interior de um conjunto.

\begin{prop}
Se $(E, \Omega)$ é um espaço pseudo-topológico, então: $A \subseteq B \Rightarrow \Upsilon(A) \subseteq \Upsilon(B)$. \\
\begin{dem}
Se $\Upsilon(A) = \emptyset$, então $\Upsilon(A) \subseteq \Upsilon(B)$. Se $\Upsilon(A)$ é um aberto D,
então $D \subseteq A \subseteq B$. Logo, $\Upsilon(A) \subseteq \Upsilon(B)$, pois $\Upsilon(B)$ é um aberto maior que D ou o próprio D.
\end{dem}
\end{prop}

\begin{prop}
Se $(E, \Omega)$ é um espaço pseudo-topológico, então: $A \subseteq B \Rightarrow \Upsilon(A) \subseteq \Upsilon(B)$ é equivalente à $\Upsilon(A) \subseteq \Upsilon(A \cup B)$. \\
\begin{dem}
$(\Rightarrow)$ Como $A \subseteq A \cup B$, então, pela hipótese, $\Upsilon(A) \subseteq \Upsilon(A \cup B)$. \\
$(\Leftarrow)$ Se $A \subseteq B$, então $A \cup B = B$. Daí, $\Upsilon(A) \subseteq \Upsilon(A \cup B) = \Upsilon(B)$.
\end{dem}
\end{prop}

\begin{prop}
Se $(E, \Omega)$ é um espaço pseudo-topológico, então: $\Upsilon(A) \cap \Upsilon(B) \subseteq \Upsilon(A \cap B)$. \\
\begin{dem}
Se $\Upsilon(A) = \emptyset$ ou $\Upsilon(B) = \emptyset$, então $\emptyset = \Upsilon(A) \cap \Upsilon(B) \subseteq \Upsilon(A \cap B)$. Agora, se $\Upsilon(A) \neq \emptyset$ e $\Upsilon(B) \neq \emptyset$, $\Upsilon(A)$ e $\Upsilon(B)$ são abertos. Logo, $\Upsilon(A) \cap \Upsilon(B)$ é um aberto e $\Upsilon(A) \cap \Upsilon(B) \subseteq A \cap B$. Por outro lado $\Upsilon(A \cap B)$ é o maior aberto contido em $A \cap B$. Assim, $\Upsilon(A) \cap \Upsilon(B) \subseteq \Upsilon(A \cap B)$.
\end{dem}
\end{prop}

Faremos uma pequena apresentação para deixarmos claro qual é a estrutura que usaremos para a interpretação da lógica do plausível.

\begin{defn}
Seja $\mathsf{A}$ uma estrutura clássica de primeira ordem com universo A. Dizemos que uma estrutura pseudo-topológica $\mathsf{K}$ para a lógica do plausível $\mathcal{L}(U)$ consiste da estrutura usual de primeira ordem, $\mathsf{A}$, acrescida de uma pseudo-topologia $(A, \Omega)$, sobre o domínio de $\mathsf{A}$.
\end{defn}

\begin{defn}
Na estrutura $\mathsf{K}$, a relação de satisfação de $\mathcal{L}(U)$ é definida de modo recursivo, no caminho usualmente utilizado, adicionando as seguintes cláusulas: \\

- Seja A uma fórmula cujas variáveis livres estão contidas em: \\

$\lbrace x \rbrace$ $\cup$ $\lbrace y_1, y_2, y_3, ... , y_n \rbrace$ e $a = (a_1, a_2, a_3, ..., a_n)$ uma sequência de elementos de A. \\

Então $\mathsf{K} \vDash Ux A[x, a] \Leftrightarrow \lbrace b \in A$ : $\mathsf{K} \vDash A[b, a] \rbrace \in \Omega$. Como usualmente, para a sentença UxAx, temos: $\mathsf{K} \vDash UxAx \Leftrightarrow \lbrace a \in A$ : $\mathsf{K} \vDash A(a) \rbrace \in \Omega$.
\end{defn}

Outras noções semânticas como modelos, validade, consequência semântica, entre outras, para $\mathcal{L}(U)$, são definidas de modo semelhante à lógica clássica de primeira ordem.

Utilizando-se desta estrutura, $\cite{gec08}$ mostraram que a lógica $\mathcal{L}(U)$ é correta e completa com respeito às estruturas pseudo-topológicas.


\section{Lógica Algébrica}

${}$ \hspace{0,5 cm} Nesta seção apresentaremos alguns conceitos sobre a lógica algébrica,
segundo as motivações de Paul Richard Halmos, mediante os textos de $\cite{hali56}$, $\cite{hal62}$ e $\cite{heg98}$.

No nosso caso, a motivação de tal seção se dá pela apresentação de Halmos de uma interpretação em modelo algébrico dos quantificadores clássicos, $\forall$ e $\exists$. 

Posteriormente, apresentaremos uma versão algébrica para o quantificador da ubiquidade, apresentado no capítulo anterior.

Para a apresentação de resultados referentes à este capítulo e, ao capítulo
posterior, utilizam-se os conceitos de ideais, filtros e homomorfismos booleanos.


\subsection{Funções Proposicionais}

${}$\hspace{0,5cm} Nesta subseção, apresentaremos as funções proposicionais em que a interpretação algébrica de uma proposição é um elemento de uma álgebra de Boole.

Abaixo, a definição formal de funções proposicionais.

\begin{defn}
Seja X um conjunto não vazio, denominado domínio, e $\textbf{B}$ uma álgebra de Boole. Consideremos o conjunto $\textbf{B}^X$ de todas as funções de X em $\textbf{B}$, e as seguintes operações. Se $p, q \in \textbf{B}^X$ e $x \in X$, então: \\

(i) O supremo de p e q, denotado por $p \lor q$, é definido por: $(p \lor q)(x) = p(x) \lor q(x)$; \\

(ii) O ínfimo de p e q, denotado por $p \land q$, é definido por: $(p \land q)(x) = p(x) \land q(x)$; \\

(iii) O complemento de p, denotado por $p'$, é definido por: $p'(x) = (p(x))'$. \\

(iv) O zero e a unidade de $\textbf{B}^X$ são, respectivamente, as funções constante iguais a 0 e a 1.
\end{defn}

O resultado abaixo é de suma importância e pode ser encontrado em $\cite{heg98}$.

\begin{prop}
$\textbf{P} = (\textbf{B}^X, 0, 1, \land, \lor, ')$ é uma Álgebra de Boole.
\end{prop}

E, abaixo, mais algumas definições.

\begin{defn} 
Uma ordem parcial natural de $\textbf{B}^X$, dada por: $p \leq q$ se, e somente se, $p \land q = p$. Claramente, $p \leq q$, em $\textbf{B}^X$, é justamente o caso em que, para todo $x \in X$, $p(x) \leq q(x)$, em $\textbf{B}$.
\end{defn}

\begin{defn}
Se $\textbf{A} = \textbf{B}^X$, então a subálgebra constante $\textbf{A}_0$ consiste das funções finitamente valoradas p de $\textbf{A}$, isto é, $\textbf{A}_0 = \lbrace p$ : $p(x) = p(y)$, para todos $x, y \in X \rbrace$.
\end{defn}

O interessante sobre $\textbf{B}^X$ é que esta álgebra pode ser considerada mais que uma Álgebra de Boole, pois temos a possibilidade de associar cada função $p$ de $\textbf{B}^X$ à um subconjunto de $\textbf{B}$, tal como na definição abaixo.

\begin{defn}
Seja $\textbf{B}^X$ como acima, o subconjunto R(p) de $\textbf{B}$, em que R(p) = $\lbrace p(x)$ : $x \in X \rbrace$, é denominado de a imagem da função p.
\end{defn}

\begin{defn}
Se $\textbf{A} = \textbf{B}^X$, a subálgebra $\textbf{A}_1$ consiste das funções finitamente valoradas com p funções constantes em $\textbf{A}$.
\end{defn}

Nota-se que, neste caso, se $p$ é uma função constante, então $R(p)$ é unitário, de modo que $R(p)$ é um subconjunto finito de $\textbf{B}$ para todos $x, y \in X$. E mais, é fácil ver que $\textbf{A}_0 \subseteq \textbf{A}_1$.

\begin{defn}
Se $p \in \textbf{A}_1$, chamamos de supremo dos valores de p, denotado por $\lor R(p)$, ao elemento de $\textbf{B}$, dado por: $\lor R(p) = p_1 \lor p_2 \lor ... \lor p_n$, em que $p_1, p_2, ..., p_n$ são todos os diferentes valores, em $\textbf{B}$, que a função proposicional p, de $\textbf{A}_1$, assume.
\end{defn}

A fim de não causar muitas confusões simbólicas, a partir deste momento, usaremos uma notação mais neutra para $\lor R(p)$, neste caso, denotaremos por $\lor R(p) = Q_{p(x)} = Qp$, para todo $x \in X$. Vale observar que $Q_p \in \textbf{A}_0$ e, ainda mais, $Q_p \in \textbf{A}_1$. Note que $Q_{p(x)}$ é um elemento de $\textbf{B}$ e não uma função $\textbf{B}^X$.

Agora, denotamos por $Q$ a função, de $\textbf{A}_1$ em $\textbf{A}_0$, que associa cada $p$ à sua respectiva função constante $Qp$.

E, algumas propriedades desta função $Q$: \\

$(P_1)$ $Q$ é $\textit{normalizado}$, ou seja, $Q0 = 0$; \\

$(P_2)$ $Q$ é $\textit{crescente}$, ou seja, $p \leq Qp$; \\

$(P_3)$ $Q$ é $\textit{distributivo}$ sobre $\lor$, ou seja, $Q(p \lor q) = Qp \lor Qq$; \\

$(P_4)$ $Q$ é $\textit{idempotente}$, ou seja, $Q(Qp) = Qp$; \\

$(P_5)$ Relação entre o complemento e Q: $Q(Qp)' = (Qp)'$; \\

$(P_6)$ Outra relação entre o complemento e Q: $Q(p') = (Qp)'$; \\

$(P_7)$ $Q$ é $\textit{quase multiplicativo sobre} \land$: $Q(p \land Qq) = Qp \land Qq$; \\

$(P_8)$ Relação entre $Q$ e $\land$: $Q(p \land q) = Qp \land q$; \\

As demonstrações destas propriedades podem ser encontradas nos textos de Halmos, principalmente em $\cite{heg98}$.


\subsection{Álgebras Monádicas Funcionais}

${}$\hspace{0,5cm} Nesta seção, apresentaremos algumas definições e alguns resultados que podem ser encontrados nos textos descritos no início do presente capítulo. Mais especificamente, em $\cite{heg98}$.

\begin{defn}
Denomina-se álgebra monádica funcional a qualquer subálgebra, $\textbf{C}$ de $\textbf{B}^X$, tal que para toda função $p \in \textbf{C}$ o supremo, $\lor R(p)$, e o ínfimo, $\land R(p)$, existem em $\textbf{C}$. Então, $\exists p$ e $\forall p$ são definidos por: $\exists p(x) = \lor R(p)$ e $\forall p(x) = \land R(p)$.
\end{defn}

A definição de $\land R(p)$ é análoga à definição dada para $\lor R(p)$, no item anterior, apenas trocando as ocorrências de $\lor$ por $\land$.

\begin{defn}
Os elementos $\exists p$ e $\forall p$ são denominados de quantificador existencial funcional e quantificador universal funcional.
\end{defn}

Temos as seguintes observações a serem feitas e que nos ajudam a justificar a denominação dada às definições acima.

Seja $X$ um conjunto não vazio e $\textbf{B}$ uma álgebra Booleana de todos os subconjuntos de $X$, e consideremos a álgebra funcional $\textbf{B}^X$. Nesta situação se $x \in X$ e $p \in \textbf{B}^X$, então $p(x) \subseteq X$.

Deste modo, a álgebra Booleana $\textbf{B}^X$ é naturalmente isomorfa a todos os subconjuntos do conjunto $X \times X$, conjunto dos pares ordenados dos elementos de $X$. Este isomorfismo é dado por uma função que atribui para cada $p \in \textbf{B}^X$ o conjunto $Y$, que é: $Y = \lbrace (x, y)$ : $y \in p(x) \rbrace$.

A noção intuitiva deste conjunto $Y$ corresponde à proposição: $\textit{y depende de p(x)}$.

Neste conjunto, a imagem, $R(p)$, de um elemento $p$ em $\textbf{B}^X$, é uma coleção de subconjuntos de $X$ e, assim, o supremo desta imagem, em $\textbf{B}$, é a união desta coleção, tal qual aquela definida na Teoria de Conjuntos. Disto segue que, cada valor de $\exists p$ corresponde à $\textit{existe um x tal que y pertence à p(x)}$.

Mais especificamente, podemos dizer que $Qp$ é a função constante em que para cada $x \in X$ atribui-se a união dos conjuntos da imagem de $p$, $\lor R(p)$. Desta maneira o conjunto $Qp$, no isomorfismo descrito acima, é: $Qp = \lbrace (x, y)$ : $\textit{existe um z tal que} y \in p(z) \rbrace$, isto é, este é o conjunto dos pares, de elementos em $X$, em que a primeira coordenada é um elemento $x$, e a segunda coordenada do par é dada por $\lor R(p)$.

Agora, um quantificador, neste caso existencial, pode ser definido como se segue.

\begin{defn}
Um quantificador existencial funcional, $\exists$, é uma função de uma álgebra Booleana em si mesma que satisfaz as seguintes propriedades: \\

(i) $\exists$ é normalizado, i.e., $\exists 0 = 0$; \\

(ii) $\exists$ é crescente, i.e., $p \leq \exists p$; \\

(iii) $\exists$ é quase multiplicativo $\land$, i.e., $\exists (p \land \exists q) = \exists p \land \exists q$.
\end{defn}

A definição de $\forall$ é similar às propriedades duais daquelas descritas para $\exists$.

\begin{defn}
Um quantificador universal funcional, $\forall$, é uma função de uma álgebra Booleana em si mesma que satisfaz as seguintes propriedades: \\

(i) $\forall 1 = 1$; \\

(ii) $\forall p \leq q$; \\

(iii) $\forall (p \lor \forall q) = \forall p \lor \forall q$.
\end{defn}

Abaixo, enunciaremos alguns resultados destes quantificadores.

\begin{teo}
Para o quantificador existencial funcional $\exists$, vale o seguinte: $\exists 1 = 1$.
\end{teo}

\begin{lema}
O quantificador existencial funcional, $\exists$, é idempotente.
\end{lema}

\begin{teo}
Se $p \leq \exists q$, então $\exists p \leq \exists q$.
\end{teo}

\begin{teo}
O quantificador existencial funcional, $\exists$, é monótono, isto é, se $p \leq q$, então $\exists p \leq \exists q$.
\end{teo}

\begin{teo}
Para o quantificador existencial funcional, $\exists$, vale o seguinte: $\exists(\exists p)' = (\exists p)'$.
\end{teo}

\begin{lema}
Uma condição necessária e suficiente para que um elemento p de $\textbf{A} = \textbf{B}^X$ pertença à imagem de $\exists$, isto é, que $p \in \exists(A)$, é que $\exists p = p$.
\end{lema}

\begin{teo}
A imagem de $\exists$, $\exists(A)$, é uma sub-álgebra Booleana de $\textbf{A} = \textbf{B}^X$.
\end{teo}

\begin{teo}
O quantificador existencial funcional, $\exists$, é disjuntivo, ou seja, $\exists(p \lor q) = \exists p \lor \exists q$.
\end{teo}

Tomemos aqui as seguintes operações Booleanas $p - q = p \land q'$ e $p + q = (p - q) \lor (q - p)$. Os resultados abaixo e suas respectivas provas são simples.

\begin{teo}
Para o quantificador existencial funcional, $\exists$, valem as seguintes propriedades: \\
(i) $\exists p - \exists q \leq \exists (p - q)$; \\
(ii) $\exists p + \exists q \leq \exists(p + q)$.
\end{teo}

Quando, para um quantificador vale o item (ii) deste teorema, dizemos o quantificador é $\textit{aditivo}$.

\begin{defn} 
Um operador de fecho topológico em uma álgebra Booleana $\textbf{A}$ é uma função Q, de $\textbf{A}$ em $\textbf{A}$, que satisfaz as seguintes condições: \\

(i) Q é normalizado; \\

(ii) Q é crescente; \\

(iii) Q é idempotente; \\

(iv) Q é aditivo.
\end{defn}

Como pudemos observar, o operador $\exists$ satisfaz todas as condições descritas acima. Assim, $\exists$ é um operador de fecho topológico.

Agora, mais alguns resultados e definições importantes serão enunciados abaixo.

\begin{teo}
Se $\exists$ é um operador de fecho topológico numa Álgebra Booleana $\textbf{A}$, então as seguintes condições são equivalentes: \\

(i) $\exists$ é um quantificador; \\

(ii) A imagem de $\exists$, $\exists(\textbf{A})$, é uma sub-álgebra Booleana de $\textbf{A}$; \\

(iii) $\exists(\exists p)' = (\exists p)'$, para todo $p \in \textbf{A}$.
\end{teo}

\begin{defn}
$\textbf{A}$ uma álgebra Booleana e $\textbf{B}$ uma subgebra Booleana de $\textbf{A}$. A álgebra $\textbf{B}$ é relativamente completa se, para todo $p \in \textbf{A}$, o conjunto B(p), definido por $B(p) = \lbrace q \in B$ : $p \leq q \rbrace$, possui um menor elemento, neste caso, um ínfimo.
\end{defn}

\begin{teo}
Se $\exists$ é um quantificador numa álgebra Booleana $\textbf{A}$ e $\textbf{B}$ é a imagem de $\exists$, então $\textbf{B}$ é uma sub-álgebra relativamente completa de $\textbf{A}$. E mais, se $B(p) = \lbrace q \in B$ : $p \leq q \rbrace$, então $\exists p = \land B(p)$, para todo $p \in \textbf{A}$.
\end{teo}

\begin{teo}
Se $\textbf{B}$ é uma sub-álgebra relativamente completa de uma álgebra Booleana $\textbf{A}$, então existe um único quantificador em $\textbf{A}$, com imagem $\textbf{B}$.
\end{teo}


\subsection{Álgebras Monádicas}

${}$\hspace{0,5cm} Nesta seção, apresentaremos as noções de álgebra monádica que, no nosso caso, é o item que iremos nos aprofundar e apresentar, posteriormente, uma álgebra monádica da ubiquidade, ou seja, uma álgebra monádica estendida por um operador que busca capturar as noções do quantificador da ubiquidade.

\begin{defn}
Uma álgebra monádica é uma álgebra Booleana $\textbf{A}$ acrescida de um quantificador funcional $\exists$ em $\textbf{A}$.
\end{defn}

Aqui, precisamos seguir, minimamente, os mesmos passos e definições dadas acima, ou seja, apresentaremos as definições para as álgebras monádicas dos conceitos que discorremos acima.

\begin{defn}
Um subconjunto $\textbf{B}$ de uma álgebra monádica $\textbf{A}$ é uma subálgebra monádica de $\textbf{A}$ se este determina uma sub-álgebra Booleana de $\textbf{A}$; e é uma álgebra monádica com relação ao quantificador de $\textbf{A}$.
\end{defn}

Em outras palavras, $\textbf{B}$ é uma sub-álgebra monádica de $\textbf{A}$ se, e somente se, $\exists p \in \textbf{B}$, sempre que $p \in \textbf{B}$.

\begin{defn}
Um homomorfismo monádico é uma função f de uma álgebra monádica em outra, tal que f é um homomorfismo Booleano e $f(\exists p) = \exists(f(p))$, para todo p.
\end{defn}

\begin{defn}
O núcleo (kernel) de um homomorfismo monádico é definido por $ker(f) = \lbrace p$ : $f(p) = 0 \rbrace$.
\end{defn}

O $\textit{kernel}$ de um homomorfismo monádico é também chamado de $\textit{ideal monádico}$.

As definições acima nos dizem que o núcleo, $ker(f)$, é um ideal Booleano em $\textbf{A}$, tal que $\exists p \in ker(f)$, sempre que $p \in ker(f)$. Similarmente, um filtro, $\textbf{F}$ em $\textbf{A}$ é um filtro Booleano em $\textbf{A}$, tal que $\forall p \in \textbf{F}$, sempre que $p \in \textbf{F}$.

\begin{defn}
Uma relação de congruência monádica, $\equiv$ em $\textbf{A}$ é uma relação de congruência Booleana em $\textbf{A}$, tal que: $p \equiv q \Rightarrow \exists p \equiv \exists q$.
\end{defn}

Sejam $\textbf{A}$ uma álgebra monádica e $\textit{I}$ é um ideal monádico em $\textbf{A}$, formemos a álgebra Booleana quociente, $\textbf{B} = \textbf{A}$ / $\textit{I}$, e consideremos o homomorfismo Booleano canônico $f$ de $\textbf{A}$ em $\textbf{B}$, que leva cada elemento $p$ na sua respectiva classe de equivalência, $[p]$ módulo $\textit{I}$. Segundo $\cite{heg98}$,existe um único caminho para converter $\textbf{B}$ numa álgebra monádica, de forma que $f$ seja um homomorfismo monádico com núcleo $\textit{I}$.

Para tal, definimos $\exists [p] = [\exists p]$, em que $\exists [p] \in \textbf{B}$ e $[\exists p] \in \textbf{A}$.

E, mais algumas definições.

\begin{defn}
Uma álgebra monádica $\textbf{A}$ é simples se $\lbrace0 \rbrace$ é seu único ideal próprio.
\end{defn}

\begin{defn}
Um ideal monádico I é maximal quando I é um ideal monádico próprio que não é um subconjunto próprio de qualquer outro ideal monádico próprio de $\textbf{A}$.
\end{defn}

Agora, chegamos a alguns resultados importantes.

\begin{lema}
Uma álgebra monádica é simples se, e somente se, é não trivial e seu quantificador é simples.
\end{lema}

\begin{lema}
Toda sub-álgebra de uma álgebra monádica simples é simples.
\end{lema}

A partir deste momento tomaremos $\textbf{O} = \lbrace 0, 1 \rbrace$, como a álgebra Booleana simples de dois elementos.

\begin{teo}
Uma álgebra monádica $\textbf{A}$ é simples se, e somente se, $\textbf{A}$ é isomorfa à uma álgebra monádica funcional $\textbf{O}$-valorada com um domínio não vazio.
\end{teo}

Abaixo alguns resultados sobre os ideais nas álgebras monádicas.

\begin{teo}
Seja I um ideal Booleano de $\textbf{A}$ e $I^p$ o conjunto de todos os elementos $p \in \textbf{A}$ em que $\exists p \in I$. Então, $I^*$ é um ideal monádico e $I^* \subseteq I$.
\end{teo}

\begin{teo}
Se I é um ideal maximal Booleano de $\textbf{A}$, então $I^p$ é um ideal monádico maximal.
\end{teo}

Em suma, o que os resultados acima nos dizem é que ao tomar cada ideal Booleano $I$ e $I^p$, a função do conjunto dos ideais Booleanos de $\textbf{A}$ no conjunto dos ideais monádicos de $\textbf{A}$, preserva a inclusão e leva ideais maximais em ideais monádicos maximais.

\begin{teo}
(Teorema dos Ideais Maximais para as álgebras Monádicas) Todo ideal monádico próprio numa álgebra monádica pode ser incluído em algum ideal monádico maximal.
\end{teo}

O teorema da existência de ideais monádicos maximais segue do teorema da existência de ideais maximais, assim como é realizado para as álgebras Booleanas.

\begin{teo}
(Teorema da existência para as álgebras monádicas) Se $\textbf{A}$ é uma álgebra monádica, então para todo elemento $p_0$, com $p_0 \neq $0 e $p_0 \in \textbf{A}$, existe um homomorfismo f de $\textbf{A}$ sobre uma álgebra monádica tal que $f(p_0) \neq 0$.
\end{teo}

Notemos que existe um ideal monádico maximal $I$ tal que $\exists p_0 \notin I$. Logo, segue um corolário.

\begin{coro}
Se $\textbf{A}$ é uma álgebra monádica, então para todo elemento $p_0$, com $p_0 \neq 0$ e $p_0 \in \textbf{A}$, existe um homomorfismo f de $\textbf{A}$ sobre uma álgebra monádica tal que $f(p_0) = 1$.
\end{coro}


\subsection{Lógica Monádica}

${}$ \hspace{0,5 cm} Neste item apresentaremos as lógicas monádicas com suas respectivas
definições.

\begin{defn} 
Uma lógica monádica é um par $(\textbf{A}, I)$, em que $\textbf{A}$ é uma álgebra monádica e I é um ideal monádico em $\textbf{A}$.
\end{defn}

Segundo, $\cite{heg98}$, as condições abaixo devem ser respeitadas pelas lógicas monádicas. \\

(i) se $p$ e $q$ são refutáveis, então $p \lor q$ são refutáveis; \\

(ii) se $p$ é refutável, então $p \land q$ deve ser refutável. \\

Assim, $I$ deve ser, no mínimo, um ideal em $\textbf{A}$. No entanto, isto não é suficiente, já que precisamos que a seguinte condição seja verdadeira: \\

(iii) se $p$ é refutável, então $\exists p$ é refutável. \\

A fim de solucionar tal questão seguem as definições.

\begin{defn}
Os elementos $p \in I$ são chamados de elementos refutáveis da lógica $(\textbf{A}, I)$. Se $\neg p \in I$, então p é demonstrável.
\end{defn}

\begin{defn}
Um modelo é uma lógica monádica $(\textbf{A}, I)$, em que $\textbf{A}$ é uma álgebra monádica funcional $\textbf{O}$-valorada, $\textbf{O} = \lbrace 0, 1 \rbrace$, e I o ideal trivial $\lbrace 0 \rbrace$.
\end{defn}

\begin{defn}
Uma interpretação de uma lógica monádica $(\textbf{A}, I)$ é um modelo $(\textbf{B}, \lbrace 0 \rbrace)$ em que há um homomorfismo monádico f de $\textbf{A}$ em $\textbf{B}$ para o qual $f(p) = 0$, sempre que $p \in I$.
\end{defn}

\begin{defn}
Todo elemento refutável é dito falso na interpretação. Se, um elemento $p \in \textbf{A}$ é falso em toda interpretação, então p é denominado universalmente inválido.
\end{defn}

\begin{defn}
Se todo elemento universalmente inválido é refutável, dizemos que a lógica é semanticamente completa.
\end{defn}

\begin{defn}
Dizemos que um elemento p é universalmente válido se $f(p) = 1$, para toda interpretação f, isto é, p é verdadeiro em toda interpretação.
\end{defn}


\subsection{Semisimplicidade e Adequação}
${}$ \hspace {0,5 cm} Neste último item apresentaremos o conceito de semisimplicidade.

\begin{defn}
Uma álgebra monádica $\textbf{A}$ é semisimples se a interpretação de todos os ideais maximais em $\textbf{A}$ é $\lbrace 0 \rbrace$.
\end{defn}

Agora podemos concluir que as álgebras monádicas constituem uma generalização das álgebras Booleanas. O teorema seguinte, em particular, mostra que toda álgebra Booleana é semisimples. Esta consequência é conhecida, pois é uma consequência imediata do Teorema da Representação de Stone e, este é o mais importante passo na prova do próximo teorema. A prova da presente generalização pode ser dada por uma imitação monádica de qualquer uma das usuais provas dos casos especiais Booleanos. A prova a seguir adota o procedimento alternativo de deduzir a generalização para o caso especial.

\begin{teo}
Toda álgebra monádica é semisimples.
\end{teo}

Utilizando-se exatamente destes conceitos é que Paul Halmos, em seus trabalhos, demonstra os resultados de correção e completude, usando apenas o ambiente algébrico, ou seja, da lógica monádica com relação a álgebra monádica. A saber, os resultados são os seguintes.

\begin{teo}
Uma lógica monádica $(\textbf{A}, I)$ é semanticamente correta se ela tem uma interpretação, isto é, se existe um homomorfismo $f: (\textbf{A}, I)$ $\rightarrow$ $(\textbf{B}, \lbrace 0 \rbrace)$, tal que se $p \in I$, então $f(p) = 0$.
\end{teo}

\begin{teo}
(Correção) Se p é demonstrável em $(\textbf{A}, I)$, então p é válido em $(\textbf{B}, \lbrace 0 \rbrace)$.
\end{teo}

\begin{teo}
(Completude) Toda lógica monádica é semanticamente completa.
\end{teo}

\begin{coro}
Se p é válido em $(\textbf{B}, \lbrace 0 \rbrace)$, então p é demonstrável em $(\textbf{A}, I)$.
\end{coro}

No próximo capítulo, estenderemos estes resultados para o operador da ubiquidade.


\section{Um modelo algébrico do quantificador da ubiquidade}

${}$ \hspace{0,5 cm} Neste capítulo apresentaremos o elemento original de nosso artigo, ou
seja, uma interpretação algébrica do quantificador da ubiquidade, tal qual apresentamos no capítulo anterior para os quantificadores clássicos.


\subsection{Quantificador funcional da Ubiquidade}

${}$ \hspace{0,5 cm} Neste item pretendemos apresentar um caso particular das funções proposicionais, neste caso, a função proposicional da ubiquidade que aqui denominaremos por quantificador funcional da ubiquidade.

\begin{defn}
Um quantificador funcional da ubiquidade ou operador da ubiquidade, denotado por $\Upsilon$, é uma operação de uma álgebra monádica $\textbf{M}$ em si mesma tal que: \\

(i) $\Upsilon p \land \Upsilon q \leq \Upsilon (p \land q)$; \\

(ii) $\Upsilon p \leq \Upsilon (p \lor q)$; \\

(iii) $\forall p \leq \Upsilon p \leq \exists p$.
\end{defn}

As condições acima introduzem as características do quantificador do plausível no contexto das álgebras monádicas. A expressão $\Upsilon p$ indica o interior do conjunto que interpreta p. 

Abaixo mostraremos alguns resultados deste quantificador funcional.

\begin{teo}
Para o operador da ubiquidade $\Upsilon$, vale o seguinte: $\Upsilon 1 = 1$. \\
\begin{dem}
Basta lembrarmos que $\forall 1 = 1 = \exists 1$. Daí, pela condição (iii) para $\Upsilon$ e para $p = 1$, $1 = \forall 1 \leq \Upsilon 1 \leq \exists 1 = 1$. Portanto, $\Upsilon 1 = 1$.
\end{dem}
\end{teo}

De modo semelhante, verificamos que $\Upsilon 0 = 0$. No contexto algébrico, a condição $\Upsilon 1 = 1$ indica que o universo é um conjunto aberto é verdadeiro, enquanto que $\Upsilon 0 = 0$ indica que o vazio é um conjunto aberto é falso.

\begin{teo}
Para o operador da ubiquidade $\Upsilon$, vale o seguinte: $\Upsilon p \land \Upsilon q \leq \exists (p \land q)$. \\
\begin{dem}
Da condição (i) para $\Upsilon$, obtemos que $\Upsilon p \land \Upsilon q \leq \Upsilon (p \land q)$ e pela condição (iii), obtemos que $\Upsilon (p \land q) \leq \exists (p \land q)$. Desse modo, $\Upsilon p \land \Upsilon q \leq \exists (p \land q)$.
\end{dem}
\end{teo}

\begin{lema}
Para o quantificador funcional da ubiquidade $\Upsilon$, vale o seguinte: $\Upsilon p \land \Upsilon p' = 0$. \\
\begin{dem}
Lembremos que $\exists 0 = 0$. Dos resultados anteriores obtemos $\Upsilon p \land \Upsilon p' \leq \exists (p \land p')$. Mas, como $p \land p' = 0$, então $\Upsilon p \land \Upsilon p' \leq \exists 0 = 0$. Como não temos o caso em que $\Upsilon p \land \Upsilon p' < 0$, então $\Upsilon p \land \Upsilon p' = 0$.
\end{dem}
\end{lema}

\begin{teo}
Para o operador da ubiquidade $\Upsilon$, vale o seguinte: $\Upsilon p \leq (\Upsilon p')'$. \\
\begin{dem}
Como $\textbf{M}$ é, em particular, uma álgebra de Boole, então $\Upsilon p' \land (\Upsilon p')' = 0$. Pelo Lema anterior, $\Upsilon p' \land \Upsilon p = 0$. Logo, $\Upsilon p \leq (\Upsilon p')'$.
\end{dem}
\end{teo}

\begin{teo}
Para o operador da ubiquidade $\Upsilon$, vale o seguinte: $\Upsilon (p \lor p') = 1$. \\
\begin{dem}
Basta lembrarmos que $p \lor p' = 1$. Daí, pelo Teorema anterior, temos $\Upsilon (p \lor p') = \Upsilon 1 = 1$.
\end{dem}
\end{teo}

Como instâncias do resultado acima, temos: $\Upsilon (p \land p') = 0$, $(\Upsilon 0)' = 1$ e $(\Upsilon (p \land p'))' = 1$.

\begin{teo}
Para o quantificador funcional da ubiquidade $\Upsilon$, vale o seguinte: Se $p \leq q$, então $\Upsilon p \leq \Upsilon q$. \\
\begin{dem}
Se $p \leq q$, então $q = p \lor q$. Logo, $\Upsilon q = \Upsilon (p \lor q)$, o que nos dá $\Upsilon p \leq \Upsilon (p \lor q) = \Upsilon q$.
\end{dem}
\end{teo}

\begin{teo}
Para o quantificador funcional da ubiquidade $\Upsilon$, vale o seguinte: $\Upsilon p \land \Upsilon q \leq \Upsilon p \lor \Upsilon q \leq \Upsilon (p \lor q)$. \\
\begin{dem}
Por instâncias da condição (ii) para $\Upsilon$, $\Upsilon p \leq \Upsilon (p \lor q)$ e $\Upsilon q \leq \Upsilon (p \lor q)$. Logo, $\Upsilon p \land \Upsilon q \leq \Upsilon p \lor \Upsilon q \leq \Upsilon (p \lor q)$.
\end{dem}
\end{teo}

\begin{teo}
Para o quantificador funcional da ubiquidade $\Upsilon$, vale o seguinte: $\Upsilon (p \land q) \leq \Upsilon p$. \\
\begin{dem}
Como $p \land q \leq p$, então, pelos Teoremas anteriores, segue que $\Upsilon (p \land q) \leq \Upsilon p$.
\end{dem}
\end{teo}

\begin{coro}
Para o quantificador funcional da ubiquidade $\Upsilon$, vale o seguinte: $\Upsilon (p \land q) = \Upsilon p \land \Upsilon q$. \\
\begin{dem}
Pela condição (i) para $\Upsilon$, temos que $\Upsilon p \land \Upsilon q \leq \Upsilon (p \land q)$. Sendo assim, basta-nos mostrar que $\Upsilon (p \land q) \leq \Upsilon p \land \Upsilon q$. Pelo Teorema anterior, obtemos que $\Upsilon (p \land q) \leq \Upsilon p$ e $\Upsilon (p \land q) \leq \Upsilon q$. Logo, $\Upsilon (p \land q) \leq \Upsilon p \land \Upsilon q$.
\end{dem}
\end{coro}


\subsection{Álgebra Monádica da Ubiquidade}

${}$ \hspace{0,5 cm} Neste item apresentaremos as noções e alguns resultados de uma álgebra monádica da ubiquidade, ou seja, uma álgebra monádica que é estendida por um operador para capturar as noções do quantificador da ubiquidade, no nosso caso o quantificador funcional da ubiquidade.

\begin{defn}
Uma álgebra monádica da ubiquidade, $\textbf{U}$, é uma álgebra monádica $\textbf{M}$ acrescida do operador da ubiquidade $\Upsilon$.
\end{defn}

A próxima definição apenas relembra definição do capítulo anterior.

\begin{defn}
O operador $\exists$ é simples quando $\exists 0 = 0$ e $\exists p = 1$, sempre que $p \neq 0$.
\end{defn}

Poderíamos tentar definir o operador da ubiquidade $\Upsilon$ como simples tal qual feito acima. Porém, temos que $\forall p \leq \Upsilon p \leq \exists p$, e também, para ideais maximais não teríamos informações sobre elementos maiores que $\Upsilon p$.

Seguiremos, de perto, os mesmos passos e definições dadas anteriormente, agora para esta álgebra.

\begin{defn}
Um subconjunto $\textbf{B}$ de uma álgebra monádica da ubiquidade $\textbf{U}$ é uma sub-álgebra monádica da ubiquidade de $\textbf{U}$ se: $p \in \textbf{B} \Rightarrow \Upsilon p \in \textbf{B}$.
\end{defn}

\begin{defn}
Um homomorfismo monádico da ubiquidade é uma função f de uma álgebra monádica da ubiquidade em outra, tal que f é um homomorfismo monádico e $f(\Upsilon p) = \Upsilon (fp)$, para todo p.
\end{defn}

\begin{defn}
O núcleo (kernel) de um homomorfismo monádico da ubiquidade f é definido por $ker(f) = \lbrace p$ : $f(p) = 0 \rbrace$.
\end{defn}

\begin{defn}
Um ideal monádico da ubiquidade I é um ideal monádico I tal que: $p \in I \Rightarrow \Upsilon p \in I$.
\end{defn}

\begin{teo}
Todo ideal monádico é um ideal monádico da ubiquidade. \\
\begin{dem}
Seja I um ideal monádico. Se $p \in I$, então $\exists p \in I$ e como $\Upsilon p \leq \exists p$, então $\Upsilon p \in I$.
\end{dem}
\end{teo}

Assim, o $\textit{kernel}$ de um homomorfismo monádico da ubiquidade é um ideal monádico da ubiquidade.

Neste trabalho, assim como fizemos no capítulo anterior, trataremos apenas dos ideais.

\begin{defn}
Uma relação de congruência monádica da ubiquidade, $\equiv$, em $\textbf{U}$ é uma relação de congruência monádica em $\textbf{U}$, tal que $\Upsilon p \equiv \Upsilon q$, sempre que $p \equiv q$.
\end{defn}

Assim, como dito anteriormente, para o caso clássico, sabemos que existe um único caminho de converter $\textbf{B}$ numa álgebra monádica da ubiquidade, de forma que $f$ torna-se um homomorfismo monádico da ubiquidade com $\textit{kernel} I$.

Para tal, temos que $\Upsilon[p] = [\Upsilon p]$, em que $\Upsilon p$ é um elemento de $\textbf{U}$.

\begin{teo}
O quantificador funcional $\Upsilon$ está bem definido em $\textbf{B}$. \\
\begin{dem}
Suponhamos $p_1, p_2 \in \textbf{U}$ e que $[p_1] = [p_2]$. Então, $p_1 \lor p_2$ está em I e, mais, $\Upsilon (p_1 \lor p_2)$ também está, pois I é um ideal da ubiquidade. Disto segue que, $\Upsilon p_1 \lor \Upsilon p_2 \in I$ pois $\Upsilon p_1 \lor \Upsilon p_2 \leq \Upsilon (p_1 \lor p_2)$, donde vem que $[\Upsilon p_1] = [\Upsilon p_2]$.
\end{dem}
\end{teo}

Agora, mais algumas definições.

\begin{defn}
Uma álgebra monádica da ubiquidade $\textbf{U}$ é simples se $\lbrace 0 \rbrace$ é o único ideal monádico da ubiquidade próprio de $\textbf{U}$.
\end{defn}

\begin{defn}
Um ideal monádico da ubiquidade, I, é maximal quando I é um ideal próprio que não é um subconjunto próprio de qualquer outro ideal monádico da ubiquidade próprio.
\end{defn}

Com isto garantimos que quase todos os Teoremas e Lemas apresentados para os ideais monádicos no capítulo anterior são válidos, agora, para os ideais monádicos da ubiquidade. Desse modo, apenas apresentaremos
versões modificadas dos enunciados destes resultados introduzindo os ideais monádicos da ubiquidade ao invés dos ideais monádicos.

\begin{lema}
Uma álgebra monádica da ubiquidade é simples se, e somente se, seu quantificador $\exists$ é simples. \\
\begin{dem}
$(\Rightarrow)$ Consideremos que U é simples, $p \in U$ e $p \neq 0$. Tomemos o seguinte conjunto $I = \lbrace q$ : $q \leq \exists p \rbrace$ o qual sabemos ser um ideal monádico para o qual $\exists$ é simples. Agora, como todo ideal monádico é monádico da ubiquidade, então I atende o enunciado. \\
$(\Leftarrow)$ Suponhamos que $\exists 0 = 0$ e $\exists p = 1$, sempre que $p \neq 0$, e mais, que I é um ideal monádico da ubiquidade em U. Assim, se $p \in I$, então $\exists p \in I$, e para $p \neq 0$, então $1 \in I$. Logo, I = U. Em outras palavras, todo ideal monádico da ubiquidade não trivial em U é impróprio, ou seja, U é simples.
\end{dem}
\end{lema}

\begin{lema} 
Toda sub-álgebra de uma álgebra monádica da ubiquidade simples é simples. \\
\begin{dem}
A única álgebra Booleana simples é a álgebra de dois elementos, denotada aqui por $\textbf{O}$. Como uma álgebra monádica da ubiquidade é simples se, e somente se, seu quantificador $\exists$ é simples, então $\textbf{O}^X$ é uma álgebra monádica da ubiquidade simples sempre que X é não vazio.
\end{dem}
\end{lema}

\begin{teo}
Uma álgebra monádica da ubiquidade $\textbf{U}$ é simples se, e somente se, $\textbf{U}$ é isomorfa à uma álgebra monádica funcional $\textbf{O}$-valorada. \\
\begin{dem}
$(\Rightarrow)$ Já foi mostrado que toda álgebra $\textbf{O}$-valorada funcional com um domínio não-vazio é simples; \\
$(\Leftarrow)$ Neste caso, nos utilizaremos novamente do Teorema de Stone referente à representação de álgebras Booleanas. Se $\textbf{U}$ é uma álgebra monádica da ubiquidade simples, então $\textbf{U}$ é, em particular, uma álgebra Booleana, na qual o Teorema de Stone é aplicável. Disto segue que existe: (i) um
conjunto X; (ii) uma sub-álgebra Booleana $\textbf{B}$ de $\textbf{O}^X$; (iii) um isomorfismo Booleano f de $\textbf{U}$ em $\textbf{B}$. Assim, pelos Lemas anteriores, o quantificador $\exists$ de $\textbf{U}$ é simples e, ainda pelos Lemas anteriores o quantificador $\exists$ de $\textbf{B}$ é simples, ou seja, f preserva todos os elementos e é automaticamente um isomorfismo monádico da ubiquidade entre as álgebras monádicas da ubiquidade $\textbf{U}$ e $\textbf{B}$.
\end{dem}
\end{teo}

Os resultados seguintes, sobre os ideais e as álgebras monádicas, são semelhantes aos apresentados no Capítulo anterior.

\begin{teo}
Se I é um ideal Booleano de $\textbf{U}$ e $I^*$ é o conjunto de todos os elementos $p \in \textbf{A}$ tais que $\Upsilon \in I$, então $I^*$ é um ideal monádico da ubiquidade. \\
\begin{dem}
Já tínhamos visto ser um ideal monádico. A condição adicional assegura ser um ideal monádico da ubiquidade.
\end{dem}
\end{teo}

\begin{teo}
Se I é um ideal maximal Booleano de $\textbf{U}$, então $I^*$ é um ideal maximal monádico da ubiquidade. \\
\begin{dem}
Já tínhamos visto ser um ideal maximal monádico. Como todo monádico é da ubiquidade, então vale o enunciado.
\end{dem}
\end{teo}

\begin{teo}
(Teorema dos ideais maximais) Todo ideal próprio numa álgebra monádica da ubiquidade está incluso em algum ideal maximal monádico da ubiquidade. \\
\begin{dem}
Já vimos que todo ideal monádico próprio está incluso num ideal maximal monádico. E, também sabemos que todo monádico é da ubiquidade, então vale o enunciado.
\end{dem}
\end{teo}

\begin{teo}
(Teorema da existência) Para todo elemento $p_0$, $p_0 \neq 0$ e $p_0 \in \textbf{U}$, em que $\textbf{U}$ é uma álgebra monádica da ubiquidade, existe um homomorfismo f de $\textbf{U}$ sobre uma álgebra monádica da ubiquidade simples tal que $f(p_0) \neq 0$.
\end{teo}

Portanto, podemos concluir que todo ideal monádico da ubiquidade próprio numa álgebra monádica da ubiquidade está incluso em algum ideal maximal monádico da ubiquidade, tal qual concluímos acima.


\subsection{Lógicas Monádicas da Ubiquidade}

${}$ \hspace{0,5 cm} Neste item, apresentaremos as lógicas monádicas da ubiquidade.

\begin{defn}
Uma lógica monádica da ubiquidade é um par $(\textbf{U}, I)$, em que $\textbf{U}$ é uma álgebra monádica da ubiquidade e I é um ideal monádico da ubiquidade em $\textbf{U}$.
\end{defn}

\begin{defn}
Os elementos $p \in I$ são os elementos refutáveis da lógica. E, se $p' \in I$, então p é provável.
\end{defn}

\begin{defn} 
Modelo é uma lógica monádica da ubiquidade $(\textbf{U}, I)$, em que $\textbf{U}$ é uma álgebra monádica funcional da ubiquidade $\textbf{O}$-valorada e I o ideal trivial $\lbrace 0 \rbrace$.
\end{defn}

\begin{defn}
Uma interpretação de uma lógica monádica da ubiquidade $(\textbf{U}, I)$ num modelo $(\textbf{B}, \lbrace 0 \rbrace)$ é um homomorfismo monádico f de $\textbf{U}$ em $\textbf{B}$ tal que $fp = 0$, sempre que $p \in I$.
\end{defn}

\begin{defn} 
Todo elemento refutável é dito falso na interpretação. Se, um elemento $p \in \textbf{U}$ é falso em toda interpretação, então p é denominado de universalmente inválido.
\end{defn}

\begin{defn}
Se todo elemento universalmente inválido é refutável, então a lógica é semanticamente completa.
\end{defn}

\begin{defn} Um elemento p é universalmente válido se $fp = 1$, para toda interpretação f, isto é, p é verdadeiro segundo toda interpretação.
\end{defn}


\subsection{Semisimplicidade e Adequação}

${}$ \hspace{0,5 cm} Neste último item, apresentaremos o conceito de semisimplicidade para uma álgebra monádica da ubiquidade. Bem como, ao final desta seção apresentaremos também os teoremas de correção e completude com relação à lógica monádica da ubiquidade e à álgebra monádica da ubiquidade.

\begin{defn}
Uma álgebra monádica da ubiquidade $\textbf{U}$ é semisimples se a interpretação para qualquer ideal maximal de $\textbf{U}$ é $\lbrace 0 \rbrace$.
\end{defn}

Dado um elemento qualquer $p$ de $\textbf{U}$, se ele está em algum ideal maximal, então a sua interpretação em $(\textbf{B}, \lbrace 0 \rbrace)$ tem que ser dada por um homomorfismo $f$ tal que $f(p) = 0$, isto é, os elementos refutáveis devem tomar valor 0. Isto também é equivalente a dizer que se para algum homomorfismo $f$, $f(p) \neq 0$, então existe um ideal maximal $I$ de $\textbf{U}$ tal que $p \notin I$.

\begin{teo}
(Teorema da semisimplicidade) Toda álgebra monádica da ubiquidade é semisimples. \\
\begin{dem}
Precisamos mostrar que se $\textbf{U}$ é uma álgebra monádica da ubiquidade, $p_0 \in \textbf{U}$ e $p_0 \neq 0$, então existe um ideal maximal monádico da ubiquidade I de $\textbf{U}$ tal que $p_0 \notin I$. Como $\textbf{U}$ é uma álgebra de Boole, então existe um ideal Booleano maximal $I_0$ tal que $p_0 \notin I_0$. Seja I o conjunto de todos os elementos $p \in \textbf{U}$ para os quais $\exists p \in I_0$. Então I é um ideal monádico e, portanto, monádico da ubiquidade e tal que $p_0 \notin I$. A prova da semisimplicidade precisa ser completada mostrando que I é maximal. Suponhamos que J é um ideal monádico, tal que $I \subset J$ (inclusão própria). Daí, existe $p \in J$ tal que $p \notin I$. Como J é monádico, então que $\exists p \in J$. Por outro lado, desde que $p \notin I$, então $\exists p \notin I_0$ e como $I_0$ é um ideal maximal Booleano, então $(\exists p)' \in I_0$. Daí, $\exists(\exists p)' = (\exists p)' \in I_0$ e, portanto, $(\exists p)' \in I \subset J$. Logo $1 \in J$ e, desse modo, $J = \textbf{U}$.
\end{dem}
\end{teo}

Agora, para mostrarmos os resultados de $\textit{correção}$ e $\textit{completude}$ fundamentalmente, precisamos mostrar que a álgebra da lógica do plausível é uma álgebra monádica da ubiquidade.

\begin{teo}
A álgebra de Lindenbaum da lógica do plausível é uma álgebra monádica da ubiquidade. \\
\begin{dem}
Como estamos tratando com uma álgebra de Lindenbaum de uma extensão conservativa da lógica de primeira ordem, devemos recordar que vale o seguinte: $\vdash A \rightarrow B$ $\Leftrightarrow$ $\vDash [A] \leq [B]$, a transferência de elementos dedutivos para conceitos algébricos da ordem Booleana. Já sabemos que a
álgebra de uma lógica de primeira ordem monádica é uma álgebra monádica. Precisamos verificar que valem as seguintes condições: \\

(i) $[\Upsilon p] \land [\Upsilon q] \leq [\Upsilon (p \land q)]$; \\

(ii) $[\Upsilon p] \leq [\Upsilon (p \lor q)]$; \\

(iii) $[\forall p] \leq [\Upsilon p] \leq [\exists p]$. \\

A condição (i) segue do $(Ax_1)$ $(\Upsilon xAx \land \Upsilon xBx) \rightarrow \Upsilon x(Ax \land Bx)$, a condição (iii) dos $(Ax_3)$ $\forall xAx \rightarrow \Upsilon xAx$ e $(Ax_4)$ $\Upsilon xAx \rightarrow \exists xAx$. A lógica do plausível é correta e completa para os modelos pseudo topológicos, conforme mostrou $\cite{gra99}$. Sabemos que o interior de A está contido no interior de $A \cup B$. Assim, a condição (ii) vale em todo espaço pseudo-topológico e, portanto, vale a sua versão correspondente na
lógica do plausível. Desse modo, a álgebra de Lindenbaum da lógica do plausível tem as propriedades que definem uma álgebra monádica da ubiquidade.
\end{dem}
\end{teo}

\begin{teo}
Uma lógica monádica da ubiquidade $(\textbf{U}, I)$ é semanticamente correta se ela tem uma interpretação, isto é, se existe um homomorfismo $f: (\textbf{U}, I) \rightarrow (\textbf{B}, \lbrace 0 \rbrace)$, tal que se $p \in I$, então $f(p) = 0$. \\
\begin{dem}
Análoga àquela dada para as lógicas monádicas, fazendo as permutas de $\textbf{M}$ por $\textbf{U}$ e I ideal monádico por I ideal monádico da ubiquidade.
\end{dem}
\end{teo}

\begin{coro}
Se p é demonstrável em $(\textbf{U}, I)$, então p é válido em $(\textbf{B}, \lbrace 0 \rbrace)$. \\
\begin{dem}
Se p é demonstrável em $(\textbf{U}, I)$, então $p' \in I$ e daí $f(p') = 0$ e, portanto, $f(p) = 1$.
\end{dem}
\end{coro}

\begin{teo}
(Completude) Toda lógica monádica da ubiquidade é semanticamente completa. \\
\begin{dem}
Desde que cada interpretação de $(\textbf{U}, I)$ num modelo $(\textbf{B}, \lbrace 0 \rbrace)$ induz, de forma natural, um homomorfismo de $\textbf{U} / I$ em $\textbf{B}$ e como a única restrição sobre $\textbf{B}$ é que seja uma álgebra simples, então a questão da completude semântica restringe-se em mostrar que $\textbf{U} / I$ é semisimples. Mas, sabemos que toda álgebra monádica da ubiquidade é semisimples e, dessa maneira, toda lógica monádica da ubiquidade é semanticamente completa. 
\end{dem}
\end{teo}

\begin{coro}
Se p é válido em $(\textbf{B}, \lbrace 0 \rbrace)$, então p é demonstrável em $(\textbf{U}, I)$. \\
\begin{dem}
O teorema anterior mostra que se p é refutável em $(\textbf{U}, I)$, então p não é válido em $(\textbf{B}, \lbrace 0 \rbrace)$.
\end{dem}
\end{coro}

Com isto concluímos que a lógica monádica da ubiquidade é correta e completa.

Como mostramos que toda álgebra de Lindenbaum da lógica do plausível é uma álgebra monádica da ubiquidade, então temos também a adequação da lógica do plausível segundo os modelos algébricos das álgebras monádicas da ubiquidade.


\section*{Considerações Finais}

${}$ \hspace {0,5 cm} O presente artigo mostrou um modelo algébrico do quantificador da ubiquidade, que foi apresentado inicialmente em versão modulada.

Para tal, discorremos, sucintamente, sobre a lógica da ubiquidade, que introduz o quantificador modulado da ubiquidade, objeto de estudo deste. A lógica da ubiquidade foi inicialmente apresentada por $\cite{gra99}$ e $\cite{gec08}$.

Prosseguindo, apresentamos as motivações dos modelos algébricos de Paul Richard Halmos, a saber, as álgebras e lógicas monádicas.

E por fim, estendemos tais noções para os elementos particulares da lógica da ubiquidade através de modelos algébricos no estilo de Halmos. Para tanto, trabalhamos apenas com uma versão de lógica monádica.

Contudo, obviamente não temos claro que seja o único, nem o melhor modelo algébrico para esta lógica. As ideias em modelos algébricos usam com bastante frequência as noções de filtros e ideais, no nosso artigo decidimos por utilizar os ideais por ser um caminho mais usual.

Assim, este trabalho torna possível verificar se uma fórmula A é ou não um teorema da lógica da ubiquidade, através de um modelo puramente algébrico.

Podemos dizer que um dos intuitos deste trabalho seria a obtenção de um modelo para a lógica da ubiquidade um pouco mais próximo à nossa intuição, pois à medida que tentamos capturar as noções de quantificador a partir de conceitos algébricos, tomando por exemplo, o quantificador existencial $(\exists)$ como uma generalização da disjunção $(\lor)$, ou seja, uma disjunção infinita, parece-nos então que, esta abordagem aproxima-se do que intuitivamente pensamos ser um quantificador.

Ao apresentarmos tais conceitos, mostramos a correção e completude deste modelo algébrico. 

É exatamente neste ponto que vemos a maior vantagem nesta abordagem, pois mostramos a adequação sem sequer sairmos do ambiente algébrico, o que nos permite certa facilidade na abordagem com relação aos métodos tradicionais.


\section*{Agradecimentos}

${}$ \hspace {0,5 cm} Agradecemos apoio da FAPESP e do DM da UNESP - Câmpus de Bauru.


\end{document}